\documentclass{groups13}
\usepackage[latin1]{inputenc}         
\usepackage{tikz}
\usetikzlibrary{arrows, decorations.markings, decorations.pathmorphing, backgrounds, positioning, fit, petri}
\newcommand{\praes}[2]{\mbox{$\langle x_1,\ldots,x_n\qua |\qua #1_1,\ldots,#1_#2\rangle$}}
\newcommand{\prae}[2]{\mbox{$\langle #1\qua |\qua #2\rangle $}}

\newcommand{\qua}{\enspace}           
\newcommand{\bewende}{\hspace*{\fill}$\Box $\\ \vspace{0.5cm}} 
\newcommand{\bewanf}{\noindent {\bf Proof:}\quad }    

\newcommand{\be}{\begin{equation}}
\newcommand{\ee}{\end{equation}}
\newcommand{\bl}[1]{\begin{equation}\label{#1}}
\begin{document}
\maketitle{ASPHERICAL WORD LABELED ORIENTED GRAPHS AND CYCLICALLY PRESENTED GROUPS}
{Jens Harlander$^{\ast}$ and Stephan Rosebrock$^{\dagger}$}
{$^{\ast}$Boise State University, Boise, Idaho, USA\newline
Email: jensharlander@boisestate.edu

$^{\dagger}$Pädagogische Hochschule Karlsruhe, Bismarckstr. 10, 76133 Karlsruhe, Germany\newline
Email: rosebrock@ph-karlsruhe.de\\[2pt]}

\begin{abstract}
A {\em word labeled oriented graph} (WLOG) is an oriented graph $\cal G$ on vertices $X=\{ x_1,\ldots ,x_k\}$, where
each oriented edge is labeled by a word in $X^{\pm1}$. 
WLOGs give rise to presentations which generalize Wirtinger presentations of knots. WLOG presentations, where the underlying graph is a tree
are of central importance in view of Whitehead's Asphericity Conjecture. We present a 
class of aspherical world labeled oriented graphs. This class can be used to produce highly non-injective aspherical labeled oriented trees and also aspherical cyclically presented groups.\\



\end{abstract}
\section{Introduction}

A {\em word labeled oriented graph} (WLOG) is an oriented graph $\cal G$ on vertices $X=\{ x_1,\ldots ,x_k\}$, where
each oriented edge is labeled by a non-empty word in $X^{\pm1}$. Associated with a word labeled oriented graph comes a {\em WLOG-presentation} $P({\cal G})$ on generators $X$ and relators in one-to-one correspondence with edges in $\cal G$. For an edge with initial vertex $x_i$, terminal vertex $x_j$, labeled $w$, we have a relation $x_iw=wx_j$. A {\em WLOG-complex} $K({\cal G})$ is the standard 2-complex associated with the WLOG-presentation $P({\cal G})$, and a {\em WLOG-group} $G({\cal G})$ is the group defined by the WLOG-presentation. We say a word labeled oriented graph is {\em aspherical} if its associated WLOG-complex is aspherical. A word labeled oriented graph is called {\em injective} if no edge label or its inverse is a subword of any other edge label. A word labeled oriented tree (WLOT) is a word labeled oriented graph where the underlying graph is a tree.\\

If every edge label of a word labeled oriented graph consists of a single letter from $X$, we simply speak of a labeled oriented graph (LOG).
LOG presentations generalize Wiritnger presentations of knots.
Labeled oriented trees (LOTs) are of central importance in view of Whitehead's Asphericity Conjecture which states that a subcomplex of an aspherical 2-complex is aspherical. LOT presentations which are Wiritinger presentations of knots are known to be aspherical by a classical result of Papakyriakopoulos \cite{Papa}.
See Bogley \cite{Bogley} and Rosebrock 
\cite{Rosebrock2} for surveys on the Whitehead Conjecture. Recently the authors have shown that injective labeled oriented trees are aspherical \cite{HR13}. Here we show that certain injective word labeled oriented graphs are aspherical (in fact, diagrammatically reducible, a strong combinatorial version of asphericity). We should note that a word 
labeled oriented graph $\cal G$ can be subdivided to produce a labeled 
oriented graph $\cal G'$, and $K({\cal G}')$ is a subdivision of $K({\cal G})$. See Example \ref{ex1} and \cite{Howie} for details. In particular, if $\cal G$ is aspherical 
then so is ${\cal G}'$. The labeled oriented graph $\cal G'$ is typically highly non-injective, even if one starts with an injective WLOG $\cal G$. Thus our result also gives access to wide classes of non-injective aspherical labeled oriented graphs and trees.  

Word labeled oriented trees have appeared before in work of Howie \cite{Howie}, where they were called weakly labeled oriented trees. \\

Another motivation for this paper comes from the interest in cyclically presented groups.
Let $F$ be the free group on generators $\{ x_1,\ldots ,x_n\}$. Let $\phi $ be the
automorphism on $F$ defined by $\phi(x_i)=x_{i+1}$ in case $1\le i \le n-1$, and $\phi(x_n)=x_{1}$.
According to Howie and Williams \cite{HW12} for $w\in F$
$$P(n,w)=\langle x_1,\ldots ,x_n\, |\, w,\phi(w),\ldots ,\phi^{n-1}(w)\rangle $$
is called a {\it cyclic presentation} and the group $G(n,w)$ it defines is called a {\it cyclically presented group}. The automorphism $\phi $ induces a {\it shift automorphism} of $G(n,w)$. There is a connection between asphericity
and the dynamics of the shift automorphism, see Bogley \cite{Bog13}.
We present several classes of cyclically presented groups which arise from word labeled oriented graphs and are aspherical
by our main result Theorem \ref{hauptsatz}.

\section{Main Theorem}

Let $K$ be the standard 2-complex associated with a finite group presentation. The complex $K$ is called {\em diagrammatically reducible} (DR) if there does not exist a reduced spherical diagram $f\colon C\to K$. See Gersten \cite{Ger87} for definitions and details. The property DR implies asphericity. 

Combinatorial curvature plays an important role in the study of asphericity of 2-complexes. Given a closed surface $S$ with a cell structure, one can assign real numbers $\omega(c)$ to the corners $c$ of the 2-cells of $S$, thought of as angles. The curvature $\kappa(v)$ at a vertex is defined as $2-\sum\omega(w)$, where the sum is taken over the corners at $v$. The curvature $\kappa(d)$ of a 2-cell $d$ is defined as $\sum \omega(c)-(|\partial d |-2)$, where the sum is taken over all the corners of $d$ and $|\partial d |$ denotes the number of edges in the boundary of $d$. The combinatorial Gauss-Bonnet Theorem asserts that summing up the curvature at all the vertices and 2-cells yields twice the Euler characteristic of the surface $S$. The combinatorial Gauss-Bonnet theorem is the basis for asphericity tests such as Gersten's weight test (see \cite{Ger87}). More details on combinatorial curvature can be found in Section 14 of McCammond's survey \cite{McCammond}.

The remaining of this section is devoted to the proof of our main result.

\begin{thm}\label{hauptsatz} Let $P(\cal G)=$\praes{r}{k} be a WLOG-presentation coming from an injective word labeled oriented graph $\cal G$. Then each relation $r_i$ is of the form
$x_{\alpha(i)}w_i=w_ix_{\beta(i)}$ where $w_i$ is a word $w_i=t_{i,1}\ldots t_{i,s_i}$ with
$t_{i,j}\in \{ x_1,\ldots ,x_n\}^{\pm 1}$. We assume $s_i\ge 2$ for $i=1,\ldots ,k$. We further assume that
\begin{enumerate}
\item Each relator is cyclically reduced.
\item For each relator  $r_i$ the words $x_{\alpha(i)} t_{i,1},\, t_{i,1}^{-1}x_{\alpha(i)},\, x_{\beta(i)} t_{i,s_i}^{-1},\, t_{i,s_i} x_{\beta(i)}$ are not pieces
(i.e.\ these words and their inverses are not subwords in another relator or in the same
relator at another place). 
\item No word $w_i$ has the form $x_j^{\pm m}$ for $m\ge 2$.
\end{enumerate}
Then $K(\cal G)$ is DR.
\end{thm}

\bewanf Let $K=K(\cal G)$. Assume there exists a reduced spherical diagram $f\colon C\to K$. We will
assign weights to the corners of the cell decomposition of $C$ such that the curvature at each vertex and at each 2-cell is non-positive. Thus the total curvature of $C$ is non-positive, contradicting the fact that $C$ is a sphere with total curvature $4$, according to Gauss-Bonnet.

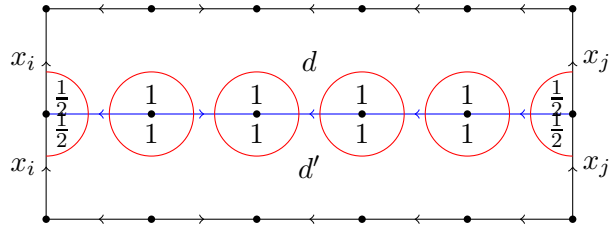
\begin{figure}[ht]\centering
\begin{tikzpicture}[scale=0.7]
\begin{scope}[decoration={markings, mark=at position 0.5 with {\arrow{>}}}]
\draw [postaction={decorate}] (0,0)--(0,2);
\draw [postaction={decorate}] (2,2)--(0,2);
\draw [postaction={decorate}] (6,2)--(4,2);
\draw [postaction={decorate}] (8,2)--(6,2);
\draw [postaction={decorate}] (10,2)--(8,2);
\draw[blue] [postaction={decorate}] (2,0)--(0,0);
\draw[blue] [postaction={decorate}] (2,0)--(4,0);
\draw[blue] [postaction={decorate}] (6,0)--(4,0);
\draw[blue] [postaction={decorate}] (8,0)--(6,0);
\draw[blue] [postaction={decorate}] (10,0)--(8,0);
\draw [postaction={decorate}] (0,-2)--(0,0);
\draw [postaction={decorate}] (2,-2)--(0,-2);
\draw [postaction={decorate}] (2,-2)--(4,-2);
\draw [postaction={decorate}] (6,-2)--(4,-2);
\draw [postaction={decorate}] (8,-2)--(6,-2);
\draw [postaction={decorate}] (10,-2)--(8,-2);
\draw [postaction={decorate}] (10,-2)--(10,0);
\draw [postaction={decorate}] (2,2)--(4,2);
\draw [postaction={decorate}] (10,0)--(10,2);
\end{scope}

\fill (0,0) circle (2pt);
\fill (0,2) circle (2pt);
\fill (2,2) circle (2pt);
\fill (4,2) circle (2pt);
\fill (6,2) circle (2pt);
\fill (8,2) circle (2pt);
\fill (10,2) circle (2pt);
\fill (2,0) circle (2pt);
\fill (4,0) circle (2pt);
\fill (6,0) circle (2pt);
\fill (8,0) circle (2pt);
\fill (10,0) circle (2pt);
\fill (0,-2) circle (2pt);
\fill (2,-2) circle (2pt);
\fill (4,-2) circle (2pt);
\fill (6,-2) circle (2pt);
\fill (8,-2) circle (2pt);
\fill (10,-2) circle (2pt);

\draw[red] (2,0) circle (0.8cm);
\draw[red] (4,0) circle (0.8cm);
\draw[red] (6,0) circle (0.8cm);
\draw[red] (8,0) circle (0.8cm);

\draw [red ] (0,-0.8)  arc (-90:90:0.8);
\draw [red ] (10,0.8)  arc (90:270:0.8);

\node at (5,1) {$d$};
\node at (5,-1) {$d'$};

\node at (0.3,0.3) {$\frac{1}{2}$};
\node at (9.7,0.3) {$\frac{1}{2}$};

\node at (0.3,-0.3) {$\frac{1}{2}$};
\node at (9.7,-0.3) {$\frac{1}{2}$};

\node[above] at (2,0) {$1$};
\node[above] at (4,0) {$1$};
\node[above] at (6,0) {$1$};
\node[above] at (8,0) {$1$};

\node[below] at (2,0) {$1$};
\node[below] at (4,0) {$1$};
\node[below] at (6,0) {$1$};
\node[below] at (8,0) {$1$};

\node[left] at (0,1) {$x_i$};
\node[right] at (10,1) {$x_j$};

\node[left] at (0,-1) {$x_i$};
\node[right] at (10,-1) {$x_j$};

\end{tikzpicture}
\caption{\label{weights1} Weights in case I. The stack line is shown in blue.}
\end{figure}

Assume $d$ is a 2-cell in $C$ with boundary word $x_iw_qx_j^{-1}w_q^{-1}$. The {\em long sides of $d$} are the sides labeled by $w_q$, the remaining two sides labeled $x_i$ and $x_j$ are called {\em short sides of $d$}. The first and last vertices of the long sides of $d$ are called {\em extremal vertices} of $d$. Note that every 2-cell in $C$ has exactly four extremal vertices. The other vertices are referred to as {\em non-extremal}. If $d_1$ and $d_2$ are 2-cells in $C$ that share a long side, that is the intersection $d_1\cap d_2$ is a long side for both $d_1$ and $d_2$, then we call $d_1\cup d_2$ a {\em stacked pair}, and the long side in the intersection a {\em stack line}. Note that because $\cal G$ is injective and $f\colon C\to K$ is reduced, if $d_1\cup d_2$ is a stacked pair, then $d_1$ and $d_2$ are mapped to the same 2-cell in $K$ under $f$,  with the same orientations. A stacked pair is shown in Figure \ref{weights1}.\\

Consider a long side $S$ of a 2-cell $d$ in $C$. We will assign weights to the corners of $d$ along $S$ so that their sum is $\le |\partial d|/2-1$.\\

{\bf Case I:} All non-extremal vertices on $S$ have valency 2.\\

In this case the long side $S$ is a stack line. Assign to the first and the last corner of $S$ 
weight $1/2$ and assign weight 1 to all other corners (see Figure \ref{weights1}). 
Note that these weights sum up to $ |\partial d|/2-1$. \\

{\bf Case II:} There is a non-extremal vertex of $S$ of valency greater or equal to 3.\\

Assign weight $2/3$ to the corners at the first and last vertex of $S$ in $d$. Assign weight $1$ to corners at vertices of valency $2$ on $S$. Let $v$  be a non-extremal vertex on $S$ of valency $\ge 3$. If $v$  is the endpoint of a stack line then assign weight $1$ to the corner in $d$ at $v$. If $v$  is not the endpoint of a stack line assign weight $2/3$ to the corner in $d$ at $v$.
\begin{figure}[ht]
\centering
\begin{tikzpicture}[scale=0.7]
\begin{scope}[decoration={markings, mark=at position 0.5 with {\arrow{>}}}]
\draw [postaction={decorate}] (0,0)--(0,2);
\draw [postaction={decorate}] (2,2)--(0,2);
\draw [postaction={decorate}] (2,2)--(4,2);
\draw [postaction={decorate}] (6,2)--(4,2);
\draw [postaction={decorate}] (8,2)--(6,2);
\draw [postaction={decorate}] (10,2)--(8,2);
\draw [postaction={decorate}] (10,0)--(10,2);
\end{scope}

\draw (0,2)--(-1,3);
\draw (10,2)--(11, 3);

\draw (4,2)--(4, 5);
\draw[blue] (6,2)--(6,5);
\draw (8,2)--(8,5);

\fill (0,2) circle (2pt);
\fill (2,2) circle (2pt);
\fill (4,2) circle (2pt);
\fill (6,2) circle (2pt);
\fill (8,2) circle (2pt);
\fill (10,2) circle (2pt);

\node at (5,1) {$d$};

\node[left] at (0,1) {$x_i$};
\node[right] at (10,1) {$x_j$};

\draw [red ] (0,1.2)  arc (-90:0:0.8);
\draw [red ] (10,1.2)  arc (270:180:0.8);

\draw [red ] (3.2,2)  arc (180:360:0.8);
\draw [red ] (5.2,2)  arc (180:360:0.8);
\draw [red ] (7.2,2)  arc (180:360:0.8);

\node at (6.2,2.2) {$v$};
\node at (5,4) {$d_1$};
\node at (7,4) {$d_2$};

\node at (0.3, 1.7) {$\frac{2}{3}$};
\node at (9.7,1.7) {$\frac{2}{3}$};

\node[above] at (5,2) {$x_k$};
\node[above] at (7,2) {$x_k$};

\node at (4,1.7) {$\frac{2}{3}$};
\node at (8,1.7) {$\frac{2}{3}$};
\node at (6,1.7) {$1$};
\end{tikzpicture}
\caption{\label{weights} Weights in case II. The stack line is shown in blue.}
\end{figure}
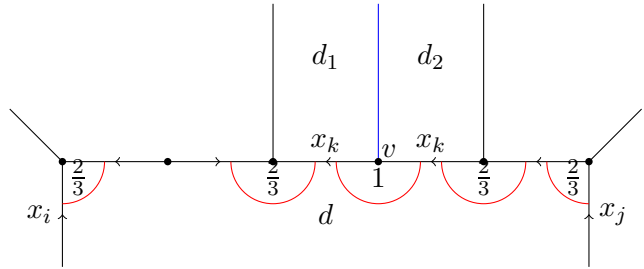

Note that there are at least 3 corners in $d$ along $S$ of weight $2/3$. There
are already two weights $2/3$ at the first and last corner of $d$ along $S$. If there were no other corners with weight $2/3$ then all non-extremal vertices along $S$ would have to be the endpoints of stack lines. But then $S$ would be labeled $x_k^{\pm m}$, for some $x_k$ and some $m\ge 2$, contradicting condition $3$ in the theorem. Thus, the sum of the weights along a long side is $\le |\partial d|/2-1$.\\

Consider a 2-cell $d$ in $C$. Using the above process we have assigned weights at the corners along both long sides of $d$. Since the sum of the weights along a long side of $d$ is  $\le |\partial d|/2-1$, the sum of the weights of all corners of $d$ is $\le |\partial d|-2$, and hence $\kappa(d)\le 0$. \\

It remains to check the curvature at the vertices of $C$. Because of condition 1 all vertices in $C$ have valency at least two. 
Because of condition 2 an extremal vertex has valency at least three.
Since the smallest weight assigned is $1/2$ we have $\kappa(v)\le 0$ for a vertex of valency greater or equal to four. Thus we only need to worry about vertices of valency two and three.
The two corners at a vertex $v$ of valency two both have weight 1 which sums up to 2. So $\kappa(v)=0$.\\

Now suppose that $v$ is a vertex of valency three. Let $d_1$, $d_2$ and $d_3$ be the three 2-cells in $C$ that share the vertex $v$ and let $c_i$ be the corner at $v$ in $d_i$.\\

If $v$ is not the end of a stack line, then the weights assigned to the corners $c_i$ are  all $2/3$. This is because weight $1/2$ is only assigned to corners at the end of stack lines and weight $1$ is assigned only to corners at the the end of stack lines or at corners at vertices of valency 2. Thus we have $\kappa(v)=0$.\\
 
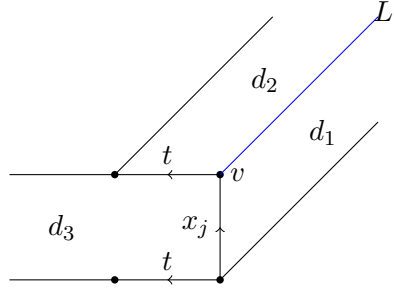
\begin{figure}[ht]
\centering
\begin{tikzpicture}[scale=0.7]
\begin{scope}[decoration={markings, mark=at position 0.5 with {\arrow{>}}}]
\draw [postaction={decorate}] (0,0)--(0,2);
\draw [postaction={decorate}] (0,2)--(-2,2);
\draw [postaction={decorate}] (0,0)--(-2,0);
\end{scope}
\draw (-2,0)--(-4,0);
\draw (-2,2)--(-4,2);
\fill (0,0) circle (2pt);
\fill (0,2) circle (2pt);
\fill (-2,2) circle (2pt);
\fill (-2,0) circle (2pt);
\draw (0,0)--(3,3);
\draw[blue] (0,2)--(3,5);
\node at (3.1,5.1) {$L$};
\draw (-2,2)--(1,5);

\node[left] at (0,1) {$x_j$};
\node[above] at (-1,0) {$t$};
\node[above] at (-1,2) {$t$};

\node at (-3,1) {$d_3$};
\node[right] at (0,2) {$v$};
\node[right] at (1.5,2.8) {$d_1$};
\node[right] at (0.4,3.8) {$d_2$};

\end{tikzpicture}
\caption{\label{exstack} Extremal vertex with stack line (blue) implies that the boundary of $d_3$ is not reduced.}
\end{figure}

Next suppose $v$ is the endpoint of a stack line $L$ and assume without loss of generality that $L=d_1\cap d_2$. Note that $d_1\cap d_3$ contains a short side of $d_1$ and $d_2\cap d_3$ contains a short side of $d_2$. Thus $L$ is the only stack line with end vertex $v$. We next argue that $v$ can not be an extremal vertex of $d_3$. Assume this would be the case. The situation is depicted in Figure \ref{exstack}. 
The stacked 2-cells $d_1$ and $d_2$ have to have the same labels on their short sides, so $t=x_j$ and the boundary word of $d$ is not reduced, contradicting the fact that each relator in $P({\cal G})$ is cyclically reduced.\\

It follows that the vertex $v$ is a non-extremal vertex of $d_3$ and hence situation is exactly as shown in Figure \ref{weights}. So the weights at $v$ are $1/2$, $1/2$, and $1$. Thus we have $\kappa(v)=0$.
\bewende

\begin{ex}\label{ex1} The WLOT-presentation 
$$P(\cal G)=\prae{a,b,c,d}{a(da)=(da)b, b(db)=(db)c, c(ac)=(ac)d}$$ satisfies the hypothesis of Theorem \ref{hauptsatz}. Hence $K(\cal G)$ is DR. Note that $\cal G$ is an interval with four vertices, each edge label has length two. Subdividing $\cal G$ into a labeled oriented interval on seven vertices, introducing a new vertex at every edge midpoint, yields a labeled oriented interval $\cal G'$ with the presentation $$P(\cal G')=\prae{a,x, b, y, c, z, d}{ad=dx, xa=ab, bd=dy, yb=bc, ca=az, zc=cd}.$$ Since $K(\cal G')$ is obtained from $K(\cal G)$ by subdivision and subdividing preserves DR (see \cite{Ger87}, 6.10), the complex $K(\cal G')$ is also DR. Note that $\cal G'$ is not injective.  \end{ex} 

\pagebreak

\section{Cyclically presented groups and further examples}

We study the cyclic presentations
$$C(n,w)=P(n,x_1wx_2^{-1}w^{-1})=$$
$$=\prae{x_1,\ldots ,x_n}{x_1w=wx_2,
x_2\phi (w)=\phi(w) x_3,\ldots ,x_n\phi^{n-1}(w)=\phi^{n-1}(w)x_1}$$
where $\phi $ is the shift automorphism. Note that $C(n,w)$ is a WLOG-presentation where the underlying graph is a circle on $n$ vertices. 

\begin{thm}\label{class} Let $w=x_j^tx_k^m$ where $t\ne 0$ and $m\ne 0$ are integers and $j\ne k$
are integers mod $n$. Assume $j\ne 1, k\ne 2, 2j-1\ne k, j\ne k-1, j\ne 2k-2, j+k\ne 3$.
Then $K(C(n,w))$ is DR.
\end{thm}

\bewanf We apply Theorem \ref{hauptsatz} to $C(n,w)$. The word labeled oriented circle that defines the WLOG-presentation $C(n,w)$ is injective because all words $w$, $\phi(w)$, $\phi^2(w)$,..., $\phi^{n-1}(w)$ are different and of the same length.
Note that the first relator $x_1x_j^tx_k^mx_2^{-1}(x_j^tx_k^m)^{-1}$ is cyclically reduced 
because of $j\ne 1, k\ne 2, j\ne k$. Since the other relators are obtained by shifting the subscripts on the first relator, all relators are cyclically reduced and hence condition 1 of the assumptions in
Theorem \ref{hauptsatz} holds. Condition 3 follows from the fact that $t \ne 0$, $m\ne 0$, and $j\ne k$.
Let $A=\{ x_1x_j^\epsilon, x_j^{-\epsilon}x_1, x_2x_k^{-\tau}, x_k^\tau x_2\}$, where $\epsilon = 1$ in case $t>0$, $\epsilon =-1$ in case $t<0$, and $\tau=1$ in case $m>0$, $\tau=-1$ in case $m<0$. Note that an element $a\in A$ being a piece would mean $a$ or $a^{-1}$ is contained in a shift of $A$, or $a$ or $a^{-1}$ is equal to a shift of $x_j^\epsilon x_k^\tau$.
The assumption $j\ne k-1$ implies $j-1\ne k-2$,
and $j+k\ne 3$ implies $k-2\ne 1-j$, which shows that $a$ or $a^{-1}$ is not a shift of an element from $A$. 
Now $j\ne 2k-2$
implies $j-k\ne k-2$, thus  $(x_2x_k^{-\tau})^{\pm 1}$ or $(x_k^\tau x_2)^{\pm 1}$ can not be a shift of $x_j^\epsilon x_k^\tau$.
Furthermore $2j-1\ne k$ implies $j-1\ne k-j$,
hence $(x_1x_j^\epsilon)^{\pm 1}$ or $(x_j^{-\epsilon}x_1)^{\pm 1}$ can not be a shift of $x_j^\epsilon x_k^\tau$. We have shown that no element from $A$ or its inverse can be a piece. Hence no shift of an element from $A$ or its inverse can be a piece and  condition 2 is indeed satisfied.  The result follows from Theorem \ref{hauptsatz}. 
\bewende

\begin{ex} The presentations $C(n,x_3^{t}x_1^{m})$, $n\ge 5$, $m\ne 0$ and $t\ne 0$,  satisfy the hypothesis of Theorem \ref{class} and are therefore DR.  One can show that all presentations $C(n,x_3x_1), n\ge 5$, do not satisfy the weight test stated in \cite{Ger87}. Thus the most obvious test for asphericity fails for these examples.
\end{ex}

Theorem \ref{class} is a result about word labeled oriented circles. Because of Whitehead's Asphericity Conjecture we are interested in
the asphericity of labeled oriented trees. We can omit one relator of $C(n,w)$ and obtain a word labeled oriented interval presentation that is also DR as a subpresentation of a presentation that is DR.

\pagebreak

\end{document}